\documentclass[11pt,twoside,reqno]{amsart}

\usepackage[a4paper,margin=2.8cm]{geometry}

\usepackage[T1]{fontenc}
\usepackage[utf8]{inputenc}
\usepackage{lmodern}
\usepackage[english]{babel}
\usepackage[final]{microtype}
\usepackage{etoolbox}

\usepackage[dvipsnames]{xcolor}
\definecolor{linkblue}{HTML}{0000FF}
\definecolor{citegreen}{HTML}{0B6B3A}
\definecolor{urlviolet}{HTML}{7A1E76}
\definecolor{linkred}{HTML}{C00000}

\usepackage{mathtools}
\usepackage{amsmath,amssymb,amsthm}
\usepackage{bbm}
\usepackage{mathrsfs}

\numberwithin{equation}{section}

\theoremstyle{plain}
\newtheorem{theorem}{Theorem}[section]
\newtheorem{proposition}[theorem]{Proposition}
\newtheorem{lemma}[theorem]{Lemma}
\newtheorem{corollary}[theorem]{Corollary}

\theoremstyle{definition}

\theoremstyle{remark}
\newtheorem{remark}[theorem]{Remark}

\usepackage{booktabs}
\usepackage{csquotes}

\usepackage[backend=biber,style=alphabetic]{biblatex}
\usepackage[
  colorlinks=true,
  linkcolor=linkred,
  citecolor=linkblue,
  urlcolor=urlviolet
]{hyperref}

\usepackage{textcase}
\makeatletter
\patchcmd{\@setauthors}{\MakeUppercase{\authors}}{\MakeTextUppercase{\authors}}{}{}
\makeatother

\usepackage[nameinlink,capitalize]{cleveref}

\crefname{theorem}{Theorem}{Theorems}
\crefname{proposition}{Proposition}{Propositions}
\crefname{lemma}{Lemma}{Lemmas}
\crefname{corollary}{Corollary}{Corollaries}
\crefname{definition}{Definition}{Definitions}
\crefname{remark}{Remark}{Remarks}

\newcommand{\R}{\mathbb{R}}
\newcommand{\T}{\mathbb{T}}
\newcommand{\eps}{\varepsilon}
\newcommand{\perpgrad}{\nabla^{\!\perp}}
\newcommand{\cof}{\operatorname{cof}}
\newcommand{\Ham}{\mathcal H}
\newcommand{\mean}[1]{\langle #1\rangle}

\newcommand{\norm}[1]{\lVert #1\rVert}

\newcommand{\citep}{\parencite}
\newcommand{\citet}{\textcite}
\DeclareMathOperator{\divg}{div}

\allowdisplaybreaks

\title[Perturbative SG--Euler limit]{%
The Semigeostrophic--Euler Limit via Perturbative Monge--Amp\`ere Estimates%
}

\author[V.~Armegioiu]{%
Victor Armegioiu\\
\NoCaseChange{Department of Mathematics, ETH Zürich}\\
\NoCaseChange{\texttt{\href{mailto:victor.armegioiu@math.ethz.ch}{victor.armegioiu@math.ethz.ch}}}%
}

\hypersetup{
  pdftitle={The Semigeostrophic-Euler Limit via Perturbative Monge-Ampere Theory},
  pdfauthor={Victor Armegioiu},
  pdfsubject={Perturbative Monge-Ampere estimates for the semigeostrophic-Euler limit},
  pdfkeywords={semigeostrophic equations, Euler limit, Monge-Ampere equation, optimal transport, Hamiltonian structure, perturbative estimates, logarithmic lifespan, stability}
}

\begin{document}

\begin{abstract}
We study the two-dimensional semigeostrophic system on the flat torus in the small-amplitude regime. We formulate the rescaled dynamics as the Lie--Poisson flow of a renormalized optimal-transport energy and expand this Hamiltonian in \(C^1\). The leading term is the Euler Hamiltonian, while the first correction is an explicit cubic Monge--Amp\`ere functional. We then derive quantitative consequences for the semigeostrophic--Euler limit: a perturbative scale-uniform endpoint Monge--Amp\`ere estimate under Hessian pinching, an explicit logarithmic perturbative lifespan for the strong branch, fixed-slow-time \(O(\eps)\) velocity convergence for canonically prepared strong branches, a conditional weak--strong rate-transfer corollary, and an \(O(\eps^2)\) Wasserstein comparison for the physical densities.
\end{abstract}

\maketitle

\begingroup
  \hypersetup{linkcolor=blue}
  \tableofcontents
\endgroup
\clearpage

\section{Introduction}
\label{sec:intro}

The semigeostrophic equations form a classical model for large-scale rotating fluids. Their dual formulation, introduced through Brenier's polar factorization and developed by Benamou--Brenier, Cullen--Gangbo, Cullen--Feldman, Loeper and others, couples incompressible transport to a Monge--Amp\`ere constraint \citep{brenier1991polar,benamou1998weak,cullen2001variational,cullen2006lagrangian,loeper2006fully,figalli2018global}. In two dimensions on the flat torus, the physical-time dual system has the form
\begin{equation*}
  \partial_s m+\nabla\cdot(m\perpgrad\psi_{\rm phys})=0,
  \qquad
  \det(I+D^2\psi_{\rm phys})=m.
\end{equation*}
Here \(m\) is the physical dual density and \(x\mapsto x+\nabla\psi_{\rm phys}(x)\) is the optimal map associated with the Monge--Amp\`ere equation.

The small-amplitude SG--Euler limit is obtained by writing
\[
  m=1+\eps\rho,
  \qquad
  \psi_{\rm phys}=\eps\psi,
  \qquad
  t=\eps s.
\]
In slow time \(t\), the equation becomes
\begin{equation*}
  \partial_t\rho+\perpgrad\psi\cdot\nabla\rho=0,
  \qquad
  \det(I+\eps D^2\psi)=1+\eps\rho.
\end{equation*}
Equivalently, in two dimensions,
\begin{equation*}
  \Delta\psi+\eps\det D^2\psi=\rho.
\end{equation*}
At leading order this is the vorticity formulation of incompressible Euler,
\begin{equation*}
  \partial_t\bar\rho+\perpgrad\bar\phi\cdot\nabla\bar\rho=0,
  \qquad
  \Delta\bar\phi=\bar\rho.
\end{equation*}

The point of this paper is to use the optimal-transport Hamiltonian structure as an organizing principle for quantitative estimates. The limit is not only a formal linearization of Monge--Amp\`ere into Poisson, and three separate quantitative questions have to be kept apart.

First, Loeper's weak modulated-energy theorem gives, in the present scaling, an estimate of the schematic form
\begin{equation}\label{eq:intro-loeper-weak}
  \mathcal E_\eps(t)\le (\mathcal E_\eps(0)+C\eps^{2/3}(1+t))e^{Ct},
\end{equation}
where \(\mathcal E_\eps\) is a squared velocity-type error \citep[Theorem~6.1]{loeper2006fully}. Thus even if the initial modulated energy is prepared at order \(O(\eps^2)\), the residual term alone yields only an \(O(\eps^{1/3})\) velocity scale. This is the correct general weak statement, but it does not exploit preparation at the level of the semigeostrophic polar factor.

Second, Loeper's strong prepared-data theorem proves fixed-slow-time convergence: for every fixed \(T>0\), \(\eps\) small enough gives a strong solution on \([0,T]\) close to the Euler branch \citep[Theorem~6.3]{loeper2006fully}. This is not the same as a quantitative lower bound on the largest admissible slow time \(T_\eps\) for a given \(\eps\). The computations used in the proof of that theorem differentiate the Monge--Amp\`ere equation along a continuity path and then apply Schauder theory to the resulting linearized equation. In those exact computations, the coefficient matrix is the cofactor of the Hessian of the physical potential. Schauder constants depend not only on ellipticity, but also on the H\"older norm of this cofactor matrix; that norm contains the H\"older seminorm of the Hessian one is trying to bound. Thus this differentiated Schauder step does not by itself provide a closed scale-uniform estimate at the quantitative level needed here. This is not a blanket objection to the fixed-time convergence theorem; it is the rate-sensitive closure issue which appears when one asks for an explicit \(\eps\)-dependent perturbative time. The perturbative endpoint estimate in \cref{sec:lifespan} supplies the anchored, scale-uniform replacement needed for that purpose.

Third, the initial data must be prepared in the variables in which SG is actually posed. One should not impose \(\psi_0^\eps=\bar\phi_0\): this generally violates the nonlinear Monge--Amp\`ere constraint. The canonical preparation is instead
\begin{equation*}
  m_0^\eps=1+\eps\bar\rho_0,
  \qquad
  \det D^2P_0^\eps=m_0^\eps,
  \qquad
  P_0^\eps=\frac{|x|^2}{2}+\eps\psi_0^\eps.
\end{equation*}
Then
\begin{equation*}
  \psi_0^\eps=\bar\phi_0+
  \eps\chi_0+O(\eps^2),
  \qquad
  \Delta\chi_0=-\det D^2\bar\phi_0,
\end{equation*}
so \(\nabla\psi_0^\eps=\nabla\bar\phi_0+O(\eps)\). For weak--strong transfer, the natural distance is even more basic: it is the physical polar-factor distance
\[
  \norm{\nabla\widetilde P_0^\eps-\nabla P_0^\eps}_{L^2}.
\]
Because \(\nabla P=x+\eps\nabla\psi\), an \(O(\eps^2)\) physical-potential mismatch is exactly an \(O(\eps)\) mismatch at the rescaled velocity level.

Our contributions are organized around these points.

\smallskip
\noindent\textbf{Hamiltonian structure.}
Set
\[
  m_\eps[\rho]=1+\eps\rho,
  \qquad
  P_\eps[\rho](x)=\frac{|x|^2}{2}+\eps\psi_\eps[\rho](x),
  \qquad
  \det D^2P_\eps=m_\eps.
\]
The map \(\nabla P_\eps\) pushes \(m_\eps dx\) onto Lebesgue measure. Hence
\[
  W_2^2(m_\eps dx,dx)
  =\int_{\T^2}|\nabla P_\eps(x)-x|^2m_\eps(x)\,dx
  =\eps^2\int_{\T^2}m_\eps|\nabla\psi_\eps|^2\,dx.
\]
The renormalized transport energy
\begin{equation*}
  \Ham_\eps(\rho):=\frac{1}{2\eps^2}W_2^2((1+\eps\rho)dx,dx)
\end{equation*}
therefore equals
\[
  \Ham_\eps(\rho)=\frac12\int_{\T^2}(1+\eps\rho)|\nabla\psi_\eps|^2\,dx.
\]
We prove that
\[
  \frac{\delta\Ham_\eps}{\delta\rho}=-\psi_\eps.
\]
With the Lie--Poisson bracket
\[
  \{F,G\}_{LP}(\rho)=\int_{\T^2}\rho\left\{\frac{\delta F}{\delta\rho},\frac{\delta G}{\delta\rho}\right\}\,dx,
  \qquad
  \{f,g\}=\perpgrad f\cdot\nabla g,
\]
the Hamiltonian equation generated by \(\Ham_\eps\) is precisely \(\partial_t\rho+\{\psi_\eps,\rho\}=0\). Thus SG is not merely a nonlinear elliptic perturbation of Euler; it is the Hamiltonian flow of a canonical optimal-transport energy.

\smallskip
\noindent\textbf{Hamiltonian expansion.}
If \(\Delta\phi=\rho\) and \(\Delta\chi=-\det D^2\phi\), then
\[
  \psi_\eps=\phi+\eps\chi+O(\eps^2)
\]
on bounded smooth sets. Consequently
\begin{equation*}
  \Ham_\eps(\rho)=\Ham_0(\rho)+\eps\Ham_1(\rho)+O(\eps^2)
  \quad\text{in }C^1,
\end{equation*}
where
\[
  \Ham_0(\rho)=\frac12\int_{\T^2}|\nabla\phi|^2\,dx,
  \qquad
  \Ham_1(\rho)=\frac13\int_{\T^2}\phi\det D^2\phi\,dx.
\]
This identifies Euler as the leading Hamiltonian limit and gives an explicit first Monge--Amp\`ere correction.

\smallskip
\noindent\textbf{Explicit lifespan beyond fixed slow time.}
The core analytic input is a perturbative endpoint estimate for the Monge--Amp\`ere equation itself. It is not a global regularity statement for arbitrary Monge--Amp\`ere solutions: it is a scale-uniform estimate on the Hessian-pinned branch. As long as \(D^2P_\eps\) remains close to the identity,
\begin{equation*}
  \norm{D^2\psi_\eps}_{L^\infty}
  \le C_\alpha \norm{\rho}_{L^\infty}
  \left(1+\log^+\frac{[\rho]_{C^\alpha}}{\norm{\rho}_{L^\infty}}\right).
\end{equation*}
This is the nonlinear Monge--Amp\`ere analogue of the endpoint logarithmic Calder\'on--Zygmund bound for Poisson. Applied to the transported perturbation, it gives
\[
  \frac{d}{dt}\norm{\nabla\rho^\eps}_{L^\infty}
  \lesssim
  \norm{\rho_0}_{L^\infty}
  \norm{\nabla\rho^\eps}_{L^\infty}
  \left(1+\log^+\frac{\norm{\nabla\rho^\eps}_{L^\infty}}{\norm{\rho_0}_{L^\infty}}\right).
\]
The transported gradient can grow double-exponentially, but the perturbative Monge--Amp\`ere condition is not \(\eps\norm{\nabla\rho^\eps}_{L^\infty}\ll1\). It is
\[
  \eps\norm{D^2\psi^\eps}_{L^\infty}\ll1.
\]
Using the endpoint estimate again, this closing condition allows
\begin{equation*}
  T_\eps\gtrsim \norm{\rho_0}_{L^\infty}^{-1}\log(1/\eps)
  \quad\text{in slow time},
  \qquad
  S_\eps=T_\eps/\eps\gtrsim \eps^{-1}\log(1/\eps)
  \quad\text{in physical time}.
\end{equation*}
In particular, this gives the weaker logarithmic-logarithmic extension one would get from the cruder gradient bootstrap, but the actual nonlinear closing gives a logarithmic slow-time lower bound. This quantifies and extends the fixed-slow-time strong convergence theory: instead of choosing \(\eps\) after a fixed \(T\), it gives an explicit lower bound for the time interval available at a given \(\eps\).

\smallskip
\noindent\textbf{Sharp prepared velocity rate and weak--strong transfer.}
For every fixed slow-time interval \([0,T]\) contained in the perturbative regime we prove
\[
  \norm{\nabla\psi^\eps(t)-\nabla\bar\phi(t)}_{L^2}\le C_T\eps
\]
for the canonically prepared strong branch. The constant is fixed-time: if one evaluates the estimate at an \(\eps\)-dependent time \(T=T_\eps\), its dependence on the Euler norms up to \(T_\eps\) has to be tracked separately. We then formulate a conditional rate-transfer result for weak or Lagrangian branches, assuming an external weak--strong stability estimate at the physical polar-factor level. If
\[
  \norm{\nabla\widetilde P^\eps(0)-\nabla P_s^\eps(0)}_{L^2}\le C\eps^2,
\]
then
\[
  \norm{\nabla\widetilde\psi^\eps(t)-\nabla\bar\phi(t)}_{L^2}\le C_T\eps.
\]
This is a conditional rate-transfer point: the general weak modulated-energy estimate loses to the residual in \eqref{eq:intro-loeper-weak}, while any weak--strong theory that supplies the physical-polar-factor estimate with constants uniform in the small-amplitude scaling preserves the strong \(O(\eps)\) velocity rate for genuinely prepared weak branches. Finally, under density pinching, the physical densities satisfy the sharper comparison
\[
  W_2(m^\eps(t),\bar m^\eps(t))\le C_T\eps^2.
\]

\section{Notation and scaling}
\label{sec:setup}

We work on \(\T^2=\R^2/\mathbb Z^2\), normalized to have unit volume. For a function \(f\), write
\[
  \mean{f}:=\int_{\T^2}f(x)\,dx,
  \qquad
  \perpgrad f=(-\partial_2f,\partial_1f),
  \qquad
  \{f,g\}=\perpgrad f\cdot\nabla g.
\]
All stream functions are normalized by zero mean. We use the convention
\begin{equation}\label{eq:euler}
  \partial_t\bar\rho+\{\bar\phi,\bar\rho\}=0,
  \qquad
  \Delta\bar\phi=\bar\rho,
  \qquad
  \mean{\bar\rho}=\mean{\bar\phi}=0.
\end{equation}
With this sign convention, the Euler velocity is \(\bar u=\perpgrad\bar\phi\).

Let \(s\) denote physical time and let \(t=\eps s\) denote slow time. We write the physical semigeostrophic density and potential as
\begin{equation*}
  m^\eps(s,x)=1+\eps\rho^\eps(t,x),
  \qquad
  P^\eps(s,x)=\frac{|x|^2}{2}+\eps\psi^\eps(t,x),
  \qquad
  t=\eps s.
\end{equation*}
The slow-time scaled SG system is
\begin{equation}\label{eq:scaled-sg}
  \partial_t\rho^\eps+\{\psi^\eps,\rho^\eps\}=0,
  \qquad
  \det(I+\eps D^2\psi^\eps)=1+\eps\rho^\eps,
  \qquad
  \mean{\rho^\eps}=\mean{\psi^\eps}=0.
\end{equation}
In two dimensions the Monge--Amp\`ere equation is equivalently
\begin{equation}\label{eq:sg-poisson-form}
  \Delta\psi^\eps+\eps\det D^2\psi^\eps=\rho^\eps.
\end{equation}
We shall repeatedly use the physical potential
\[
  P^\eps=\frac{|x|^2}{2}+\eps\psi^\eps,
  \qquad
  D^2P^\eps=I+\eps D^2\psi^\eps.
\]
The perturbative Monge--Amp\`ere regime is the set of times on which
\begin{equation}\label{eq:bootstrap}
  \norm{D^2P^\eps-I}_{L^\infty}
  =\eps\norm{D^2\psi^\eps}_{L^\infty}
  \le \eta_0,
\end{equation}
for a small universal constant \(\eta_0>0\).

For a mean-zero distribution \(f\), we use
\[
  \norm{f}_{\dot H^{-1}}:=\norm{\nabla(-\Delta)^{-1}f}_{L^2}.
\]
Notice that with our convention \(\Delta\phi=f\), one has
\[
  \norm{f}_{\dot H^{-1}}=\norm{\nabla\phi}_{L^2}.
\]

\section{Exact optimal-transport Hamiltonian structure}
\label{sec:hamiltonian}

Let \(\rho\in C^\infty(\T^2)\) be mean-zero and assume \(m_\eps=1+\eps\rho>0\). Let
\begin{equation*}
  P_\eps[\rho]=\frac{|x|^2}{2}+\eps\psi_\eps[\rho]
\end{equation*}
be the periodic optimal-transport potential satisfying
\begin{equation*}
  \det D^2P_\eps=m_\eps,
  \qquad
  \nabla P_\eps{}_{\#}(m_\eps dx)=dx,
  \qquad
  \mean{\psi_\eps}=0.
\end{equation*}
The existence and uniqueness of this potential are the periodic polar factorization theorem of Brenier--McCann; see \citep{brenier1991polar,mccann2001polar,loeper2006fully}.

We use the quadratic cost induced by the flat-torus distance \(d_{\T^2}\).  The optimal map is represented by the periodic lift
\[
  \nabla P_\eps(x)=x+\eps\nabla\psi_\eps(x),
  \qquad \nabla\psi_\eps(x+z)=\nabla\psi_\eps(x),\quad z\in\mathbb Z^2.
\]
On the perturbative branch \(\eps\norm{D^2\psi_\eps}_{L^\infty}\le\eta_0\), after decreasing \(\eta_0\) if necessary, Poincar\'e's inequality for the periodic vector field \(\eps\nabla\psi_\eps\) gives
\[
  \norm{\eps\nabla\psi_\eps}_{L^\infty}
  \le C\eps\norm{D^2\psi_\eps}_{L^\infty}
  \le C\eta_0<\frac12.
\]
Thus the lifted displacement is the minimizing torus displacement, and
\[
  d_{\T^2}(x,\nabla P_\eps(x))
  =|\eps\nabla\psi_\eps(x)|.
\]
All uses of the identity below are on this small periodic branch.  Away from this branch the same expression should be interpreted through the periodic polar-factor lift and the torus distance.

Define
\begin{equation*}
  \Ham_\eps(\rho):=\frac{1}{2\eps^2}W_2^2((1+\eps\rho)dx,dx).
\end{equation*}
Since \(\nabla P_\eps\) is the optimal map from \(m_\eps dx\) to \(dx\), and since on the perturbative lift the torus displacement is \(\eps\nabla\psi_\eps\),
\[
  W_2^2(m_\eps dx,dx)
  =\int_{\T^2}d_{\T^2}(x,\nabla P_\eps(x))^2m_\eps(x)\,dx
  =\eps^2\int_{\T^2}m_\eps|\nabla\psi_\eps|^2\,dx.
\]
Hence
\begin{equation}\label{eq:Ham-identity}
  \Ham_\eps(\rho)
  =\frac12\int_{\T^2}m_\eps|\nabla\psi_\eps|^2\,dx.
\end{equation}

\begin{theorem}[Exact first variation]\label{thm:first-variation}
For smooth mean-zero \(\rho\) with \(1+\eps\rho>0\), and with the periodic polar factor represented on the perturbative lift described above, the functional \(\Ham_\eps\) is \(C^1\) and
\begin{equation*}
  D\Ham_\eps(\rho)[h]
  =-\int_{\T^2}\psi_\eps[\rho]h\,dx
\end{equation*}
for every smooth mean-zero perturbation \(h\). Equivalently,
\[
  \frac{\delta\Ham_\eps}{\delta\rho}=-\psi_\eps[\rho].
\]
\end{theorem}

\begin{proof}
Let \(\rho_\tau=\rho+\tau h\), \(m_\tau=1+\eps\rho_\tau\), and
\[
  P_\tau=\frac{|x|^2}{2}+\eps\psi_\tau.
\]
Set
\[
  \zeta=\left.\frac{d}{d\tau}\right|_{\tau=0}\psi_\tau,
  \qquad
  M=\cof D^2P_\eps.
\]
Differentiating \(\det D^2P_\tau=m_\tau\) at \(\tau=0\) gives
\[
  M:D^2\zeta=h.
\]
By the Piola identity for Hessian cofactors,
\[
  \partial_iM_{ij}=0,
\]
and therefore
\begin{equation}\label{eq:lin-MA}
  h=\divg(M\nabla\zeta).
\end{equation}
Differentiating \eqref{eq:Ham-identity} gives
\begin{equation}\label{eq:dH-raw}
  D\Ham_\eps(\rho)[h]
  =\int_{\T^2}m_\eps\nabla\psi_\eps\cdot\nabla\zeta\,dx
  +\frac\eps2\int_{\T^2}|\nabla\psi_\eps|^2h\,dx.
\end{equation}
Using \eqref{eq:lin-MA} in the second term and integrating by parts,
\[
  \frac\eps2\int |\nabla\psi_\eps|^2h
  =-\frac\eps2\int M\nabla(|\nabla\psi_\eps|^2)\cdot\nabla\zeta.
\]
Now
\[
  D^2P_\eps\nabla\psi_\eps
  =\nabla\psi_\eps+\eps D^2\psi_\eps\nabla\psi_\eps
  =\nabla\psi_\eps+\frac\eps2\nabla(|\nabla\psi_\eps|^2).
\]
Since \(M D^2P_\eps=(\det D^2P_\eps)I=m_\eps I\), this yields
\[
  M\nabla\psi_\eps+\frac\eps2M\nabla(|\nabla\psi_\eps|^2)=m_\eps\nabla\psi_\eps.
\]
Thus \eqref{eq:dH-raw} becomes
\[
  D\Ham_\eps(\rho)[h]
  =\int M\nabla\psi_\eps\cdot\nabla\zeta
  =-\int \psi_\eps\divg(M\nabla\zeta)
  =-\int \psi_\eps h.
\]
This proves the claim.
\end{proof}

On this branch, the preceding theorem gives an exact Hamiltonian formulation. For functionals \(F,G\) on mean-zero densities define the Lie--Poisson bracket
\begin{equation}\label{eq:LP-bracket}
  \{F,G\}_{LP}(\rho)
  :=\int_{\T^2}\rho\left\{\frac{\delta F}{\delta\rho},\frac{\delta G}{\delta\rho}\right\}\,dx.
\end{equation}

\begin{corollary}[SG as a Lie--Poisson Hamiltonian flow]\label{cor:ham-flow}
The scaled semigeostrophic equation \eqref{eq:scaled-sg} is the Hamiltonian flow generated by \(\Ham_\eps\) with respect to \eqref{eq:LP-bracket}.
\end{corollary}

\begin{proof}
For the linear observable \(F_f(\rho)=\int \rho f\,dx\), \(\delta F_f/\delta\rho=f\). By \cref{thm:first-variation},
\[
  \frac{d}{dt}F_f(\rho(t))
  =\{F_f,\Ham_\eps\}_{LP}(\rho)
  =\int\rho\{f,-\psi_\eps\}
  =\int\rho\{\psi_\eps,f\}.
\]
Since \(\perpgrad\psi_\eps\) is divergence-free,
\[
  \int\rho\{\psi_\eps,f\}
  =-\int f\{\psi_\eps,\rho\}.
\]
Therefore \(\partial_t\rho+\{\psi_\eps,\rho\}=0\) in distributions.
\end{proof}

\section{Expansion of the Hamiltonian}
\label{sec:expansion}

We now expand \(\Ham_\eps\) on smooth bounded sets. Fix \(k\ge3\), \(\alpha\in(0,1)\), and work in the mean-zero spaces \(C^{k,\alpha}_0(\T^2)\). For \(\rho\in C^{k,\alpha}_0\), write \(\psi_\eps[\rho]\) for the mean-zero solution of
\begin{equation}\label{eq:psi-eps}
  \Delta\psi_\eps+\eps\det D^2\psi_\eps=\rho.
\end{equation}

\begin{lemma}[Smooth dependence and stream-function expansion]\label{lem:smooth-expansion}
Let \(B\subset C^{k,\alpha}_0(\T^2)\) be bounded. For \(|\eps|\) sufficiently small depending on \(B\), the map \(\rho\mapsto\psi_\eps[\rho]\) is smooth from \(B\) into \(C^{k+2,\alpha}_0\), and
\begin{equation*}
  \psi_\eps[\rho]=\phi[\rho]+\eps\chi[\rho]+O_B(\eps^2)
  \quad\text{in }C^{k+2,\alpha},
\end{equation*}
where
\begin{equation}\label{eq:phi-chi}
  \Delta\phi=\rho,
  \qquad
  \Delta\chi=-\det D^2\phi,
  \qquad
  \mean{\phi}=\mean{\chi}=0.
\end{equation}
The same expansion holds after one derivative in \(\rho\).
\end{lemma}

\begin{proof}
Consider
\[
  \mathcal F(\eps,\psi,\rho)=\Delta\psi+\eps\det D^2\psi-\rho.
\]
At \(\eps=0\), the derivative in \(\psi\) is \(D_\psi\mathcal F(0,\psi,\rho)=\Delta\), an isomorphism from \(C^{k+2,\alpha}_0\) to \(C^{k,\alpha}_0\). The implicit function theorem gives smooth dependence for small \(|\eps|\). Substituting \(\psi_\eps=\phi+\eps\chi+\eps^2r_\eps\) into \eqref{eq:psi-eps} and comparing powers of \(\eps\) gives \eqref{eq:phi-chi}. The equation for \(\chi\) is solvable on the torus because
\[
  \int_{\T^2}\det D^2\phi\,dx=0,
\]
which follows by writing \(\det D^2\phi\) in divergence form. The bounded remainder and the differentiated expansion follow from the same implicit-function argument.
\end{proof}

\begin{lemma}[Cubic identity]\label{lem:cubic-id}
For every smooth periodic \(\phi\),
\begin{equation}\label{eq:cubic-id}
  \int_{\T^2}\Delta\phi\,|\nabla\phi|^2\,dx
  =-\frac43\int_{\T^2}\phi\det D^2\phi\,dx.
\end{equation}
\end{lemma}

\begin{proof}
Writing subscripts for derivatives and integrating by parts,
\[
  \int\Delta\phi|\nabla\phi|^2
  =-2\int \phi_i\phi_j\phi_{ij}.
\]
On the other hand, a direct integration by parts gives
\[
  \int\phi(\phi_{11}\phi_{22}-\phi_{12}^2)
  =\frac32\int \phi_i\phi_j\phi_{ij}.
\]
Combining the two identities gives \eqref{eq:cubic-id}.
\end{proof}

\begin{theorem}[Hamiltonian expansion]\label{thm:ham-expansion}
Let \(B\subset C^{k,\alpha}_0(\T^2)\), \(k\ge3\), be bounded. Then
\begin{equation*}
  \Ham_\eps(\rho)=\Ham_0(\rho)+\eps\Ham_1(\rho)+O_B(\eps^2)
  \quad\text{in }C^1(B),
\end{equation*}
where
\begin{equation*}
  \Ham_0(\rho)=\frac12\int_{\T^2}|\nabla\phi|^2\,dx,
  \qquad
  \Ham_1(\rho)=\frac13\int_{\T^2}\phi\det D^2\phi\,dx,
  \qquad
  \Delta\phi=\rho.
\end{equation*}
Moreover
\begin{equation}\label{eq:H0-H1-vars}
  \frac{\delta\Ham_0}{\delta\rho}=-\phi,
  \qquad
  \frac{\delta\Ham_1}{\delta\rho}=-\chi,
  \qquad
  \Delta\chi=-\det D^2\phi.
\end{equation}
Consequently the SG vector field expands as
\begin{equation}\label{eq:vector-field-expansion}
  -\{\psi_\eps[\rho],\rho\}
  =-\{\phi,\rho\}-\eps\{\chi,\rho\}+O_B(\eps^2).
\end{equation}
Equivalently, the first modified equation is
\begin{equation*}
  \partial_t\rho+\{\phi,\rho\}+\eps\{\chi,\rho\}=O(\eps^2),
  \qquad
  \Delta\phi=\rho,
  \qquad
  \Delta\chi=-\det D^2\phi.
\end{equation*}
\end{theorem}

\begin{proof}
From \eqref{eq:Ham-identity},
\[
  \Ham_\eps(\rho)=\frac12\int(1+\eps\rho)|\nabla\psi_\eps|^2.
\]
Using Lemma~\ref{lem:smooth-expansion},
\[
  \psi_\eps=\phi+\eps\chi+O(\eps^2).
\]
Therefore
\[
  \Ham_\eps
  =\frac12\int|\nabla\phi|^2
  +\eps\left(\int\nabla\phi\cdot\nabla\chi+\frac12\int\rho|\nabla\phi|^2\right)
  +O(\eps^2).
\]
Since \(\Delta\chi=-\det D^2\phi\),
\[
  \int\nabla\phi\cdot\nabla\chi=-\int\phi\Delta\chi=\int\phi\det D^2\phi.
\]
Since \(\rho=\Delta\phi\), Lemma~\ref{lem:cubic-id} gives
\[
  \frac12\int\rho|\nabla\phi|^2=-\frac23\int\phi\det D^2\phi.
\]
Thus the coefficient of \(\eps\) is \(\frac13\int\phi\det D^2\phi\). The \(C^1\) statement follows either by differentiating the above expansion from Lemma~\ref{lem:smooth-expansion}, or directly from \cref{thm:first-variation}.

It remains only to verify the displayed first variations. If \(\Delta\eta=h\), then
\[
  D\Ham_0(\rho)[h]=\int\nabla\phi\cdot\nabla\eta=-\int\phi h.
\]
For \(J(\phi)=\int\phi\det D^2\phi\), repeated integration by parts gives
\[
  DJ(\phi)[\eta]=3\int\eta\det D^2\phi.
\]
Hence
\[
  D\Ham_1(\rho)[h]=\int\eta\det D^2\phi.
\]
If \(\Delta\chi=-\det D^2\phi\), then
\[
  -\int\chi h=-\int\chi\Delta\eta=-\int\eta\Delta\chi=\int\eta\det D^2\phi.
\]
This proves \eqref{eq:H0-H1-vars}. Finally \eqref{eq:vector-field-expansion} follows from \(\psi_\eps=\phi+\eps\chi+O(\eps^2)\).
\end{proof}

\section{Perturbative Monge--Amp\`ere endpoint estimate and lifespan}
\label{sec:lifespan}

The perturbative lifespan argument uses one elliptic estimate which has to be formulated carefully.  The tempting rewrite
\[
  \Delta\psi=\rho-\eps\det D^2\psi
\]
is not an endpoint estimate for the semigeostrophic Monge--Amp\`ere equation.  It would bound the Poisson right-hand side by
\[
  \norm{\rho}_{L^\infty}+
  \eps\norm{D^2\psi}_{L^\infty}^2,
\]
and therefore reintroduce a quadratic feedback in the unknown quantity.  The correct statement is an anchored estimate for the Hessian of the Monge--Amp\`ere potential itself.  In this subsection we give the full reduction, including the precise Campanato input and the large-scale anchor.

\subsection*{The cofactor--Schauder closure in prior work}

We first isolate the point in the classical strong convergence proof which motivates the replacement estimate.  The computations below are exactly the continuity-path computations used in the proof of \cite[Theorem~6.3]{loeper2006fully} to obtain the estimate labelled (47), rewritten in the present scaling.  In the present notation this path is
\begin{equation}\label{eq:loeper-path}
  \det(I+\eps D^2\psi_t)=1+t\eps\rho,
  \qquad 0\le t\le1.
\end{equation}
Differentiating \eqref{eq:loeper-path} in the path parameter gives
\begin{equation*}
  \cof(I+\eps D^2\psi_t):\eps D^2\dot\psi_t
  =\eps\rho,
\end{equation*}
where \(\dot\psi_t=\partial_t\psi_t\).  Equivalently,
\begin{equation}\label{eq:loeper-linearized-clean}
  M_t^{ij}\partial_{ij}\dot\psi_t=\rho,
  \qquad
  M_t:=\cof(I+\eps D^2\psi_t).
\end{equation}
The ellipticity of \eqref{eq:loeper-linearized-clean} is not the issue.  If
\[
  \eps\norm{D^2\psi_t}_{L^\infty}\le\eta_0,
\]
then \(M_t\) is uniformly positive definite.  The difficulty is that the Schauder estimate used for \eqref{eq:loeper-linearized-clean} also depends on the coefficient H\"older norm.  On a ball, the estimate has the schematic form
\begin{equation}\label{eq:schauder-cofactor-dependence}
  \norm{\dot\psi_t}_{C^{2,\alpha}(B_{1/2})}
  \le
  C\bigl(\lambda,\Lambda,[M_t]_{C^\alpha(B_1)}\bigr)
  \left(
    \norm{\dot\psi_t}_{L^\infty(B_1)}+
    \norm{\rho}_{C^\alpha(B_1)}
  \right).
\end{equation}
In two dimensions,
\begin{equation*}
  \cof(I+\eps D^2\psi_t)=I+
  \eps\cof(D^2\psi_t).
\end{equation*}
Explicitly, if
\[
  D^2\psi_t=
  \begin{pmatrix}
    \psi_{11} & \psi_{12}\\
    \psi_{12} & \psi_{22}
  \end{pmatrix},
\]
then
\[
  \cof(I+\eps D^2\psi_t)=
  \begin{pmatrix}
    1+\eps\psi_{22} & -\eps\psi_{12}\\
    -\eps\psi_{12} & 1+\eps\psi_{11}
  \end{pmatrix}.
\]
Consequently,
\begin{equation*}
  [M_t]_{C^\alpha}
  \le C\eps[D^2\psi_t]_{C^\alpha}.
\end{equation*}
Thus the displayed Schauder step does not by itself provide a closed scale-uniform estimate: the constant in \eqref{eq:schauder-cofactor-dependence} depends on the H\"older seminorm of the same Hessian one is trying to control.  The obstruction is not loss of ellipticity, but absence of an independent scale-uniform bound for the coefficient regularity.

The estimate below avoids differentiating the continuity path.  We set
\[
  w=\eps\psi,
\]
and use the equation in the form
\[
  \det(I+D^2w)=1+f.
\]
Equivalently,
\[
  F(D^2w)=g,
  \qquad
  F(A)=\log\det(I+A),
  \qquad
  g=\log(1+f).
\]
The Campanato part controls the oscillation of \(D^2w\) relative to quadratic polynomials.  The absolute quadratic component is then fixed by an independent \(L^2\)-anchor.  This separation is the missing point in any argument which tries to treat the cofactor Schauder estimate as if bounded ellipticity alone controlled the constants.

For a bounded function \(h\), write
\[
  \omega_h(r):=\sup_{|x-y|\le r}|h(x)-h(y)|.
\]

\begin{lemma}[Spectral concave extension]\label{lem:spectral-extension}
There exist universal constants \(\eta_0>0\), \(0<\lambda\le\Lambda<\infty\), and a smooth concave uniformly elliptic function
\[
  \widetilde F:\operatorname{Sym}(2)\to\R
\]
such that
\[
  \widetilde F(A)=\log\det(I+A)
  \qquad\text{whenever } |A|\le 2\eta_0,
\]
and for all symmetric matrices \(A\) and all positive semidefinite symmetric matrices \(B\),
\begin{equation}\label{eq:global-ellipticity-Ftilde}
  \lambda\operatorname{tr}B
  \le
  D\widetilde F(A)[B]
  \le
  \Lambda\operatorname{tr}B.
\end{equation}
\end{lemma}

\begin{proof}
Choose \(\eta_0>0\) so small that \([-2\eta_0,2\eta_0]\subset(-1/2,1/2)\).  Let
\[
  q(s)=\log(1+s)
  \qquad\text{for }s\in[-2\eta_0,2\eta_0].
\]
On this interval,
\[
  q'(s)=\frac1{1+s},
  \qquad
  q''(s)=-\frac1{(1+s)^2}<0,
\]
and therefore \(q'\) is bounded above and below by positive universal constants.  Extend \(q\) to a smooth concave function \(\tilde q:\R\to\R\) such that
\[
  \tilde q=q\quad\text{on }[-2\eta_0,2\eta_0],
  \qquad
  0<\lambda\le \tilde q'(s)\le\Lambda<\infty,
  \qquad
  \tilde q''(s)\le0
\]
for all \(s\in\R\).  This is obtained by extending the positive decreasing derivative \(q'\) from the small interval to a smooth positive decreasing function bounded between \(\lambda\) and \(\Lambda\), and then integrating.

For a symmetric matrix \(A\), let \(\lambda_1(A),\lambda_2(A)\) be its eigenvalues and define
\[
  \widetilde F(A)=\tilde q(\lambda_1(A))+\tilde q(\lambda_2(A)).
\]
This spectral function is concave because \(\tilde q\) is concave, and it is uniformly elliptic because \(\tilde q'\) is bounded above and below.  More explicitly, if \(B\ge0\), the first variation is
\[
  D\widetilde F(A)[B]
  =\operatorname{tr}(\tilde q'(A)B),
\]
where \(\tilde q'(A)\) is defined by the spectral calculus.  Since
\[
  \lambda I\le \tilde q'(A)\le \Lambda I,
\]
we get \eqref{eq:global-ellipticity-Ftilde}.  Finally, if \(|A|\le2\eta_0\), then both eigenvalues of \(A\) lie in \([-2\eta_0,2\eta_0]\), and hence
\[
  \widetilde F(A)=\sum_{i=1}^2\log(1+\lambda_i(A))=
  \log\det(I+A).
\]
\end{proof}

We now prove the local endpoint estimate used below.  The proof is included in full because this is the rate-sensitive point of the paper.  We use only the standard compactness form of the Evans--Krylov--ABP improvement step for concave uniformly elliptic equations, stated next.  The subsequent dyadic iteration, anchoring, logarithmic Dini evaluation, and perturbative Monge--Amp\`ere scaling are all written out explicitly.

\begin{lemma}[One-step quadratic improvement]\label{lem:one-step-improvement}
Fix ellipticity constants \(0<\lambda\le\Lambda\).  There exist constants
\[
  \gamma_*\in(0,1),\qquad
  \theta\in(0,1/4),\qquad
  \delta_*>0,
\]
depending only on \(\lambda,\Lambda\), with the following property.  Let
\(F:\operatorname{Sym}(2)\to\R\) be concave, uniformly elliptic with constants
\(\lambda,\Lambda\), independent of \(x\), and normalized by \(F(0)=0\).  Let
\(v\in C^2(B_1)\) solve
\[
  F(D^2v)=h
  \qquad\text{in }B_1,
  \qquad h(0)=0,
\]
and assume
\[
  \norm{v}_{L^\infty(B_1)}\le1,
  \qquad
  \norm{h}_{L^\infty(B_1)}\le\delta_*.
\]
Then there is a quadratic polynomial
\[
  P(x)=a+b\cdot x+\frac12 x^TAx,
  \qquad F(A)=0,
\]
such that
\begin{equation}\label{eq:one-step-excess}
  \norm{v-P}_{L^\infty(B_\theta)}\le \theta^{2+\gamma_*},
\end{equation}
and
\begin{equation}\label{eq:one-step-coeff}
  |a|+|b|+|A|\le C_*.
\end{equation}
Here \(C_*\) depends only on \(\lambda,\Lambda\).
\end{lemma}

\begin{proof}
We recall the standard compactness proof, since the constants have to be uniform in later rescalings.  Suppose the conclusion is false.  Then there are concave uniformly elliptic operators \(F_n\), functions \(v_n\), and right-hand sides \(h_n\) with
\[
  \norm{v_n}_{L^\infty(B_1)}\le1,
  \qquad
  F_n(D^2v_n)=h_n,
  \qquad
  h_n(0)=0,
  \qquad
  \norm{h_n}_{L^\infty(B_1)}\to0,
\]
for which no polynomial satisfying \eqref{eq:one-step-excess}--\eqref{eq:one-step-coeff} exists.

By the ABP stability and interior H\"older estimates for uniformly elliptic concave equations \citep{caffarelli1995fully}, after passing to a subsequence we have
\[
  v_n\to v_\infty\quad\text{locally uniformly in }B_1,
\]
and the operators \(F_n\), normalized by \(F_n(0)=0\), converge locally uniformly on compact subsets of \(\operatorname{Sym}(2)\) to a concave uniformly elliptic operator \(F_\infty\) with the same ellipticity constants.  The limit satisfies
\[
  F_\infty(D^2v_\infty)=0
  \qquad\text{in }B_1
\]
in the viscosity sense.  By the Evans--Krylov theorem for concave uniformly elliptic equations \citep{caffarelli1995fully}, there are \(\bar\gamma\in(0,1)\) and \(C\), depending only on \(\lambda,\Lambda\), such that
\[
  \norm{v_\infty}_{C^{2,\bar\gamma}(B_{3/4})}\le C.
\]
Choose \(\gamma_*\in(0,\bar\gamma)\) and then choose \(\theta\in(0,1/4)\) so small that
\[
  C\theta^{2+\bar\gamma}\le \frac12\theta^{2+\gamma_*}.
\]
Let \(P_\infty\) be the second-order Taylor polynomial of \(v_\infty\) at the origin:
\[
  P_\infty(x)=v_\infty(0)+\nabla v_\infty(0)\cdot x
  +\frac12 x^TD^2v_\infty(0)x.
\]
Then
\[
  F_\infty(D^2P_\infty)=F_\infty(D^2v_\infty(0))=0,
\]
\[
  \norm{v_\infty-P_\infty}_{L^\infty(B_\theta)}
  \le \frac12\theta^{2+\gamma_*},
\]
and the coefficients of \(P_\infty\) are bounded by a constant depending only on \(\lambda,\Lambda\).

For large \(n\), local uniform convergence gives
\[
  \norm{v_n-P_\infty}_{L^\infty(B_\theta)}
  \le \theta^{2+\gamma_*}.
\]
The only remaining point is that the frozen polynomial must satisfy \(F_n(D^2P_n)=0\), not merely \(F_\infty(D^2P_\infty)=0\).  Because \(F_n\to F_\infty\) locally uniformly and each \(F_n\) is uniformly elliptic, the scalar map
\[
  t\mapsto F_n(D^2P_\infty+tI)
\]
is strictly increasing with slope between \(2\lambda\) and \(2\Lambda\).  Therefore there is a unique \(t_n\to0\) such that
\[
  F_n(D^2P_\infty+t_nI)=0.
\]
Set
\[
  P_n(x)=P_\infty(x)+\frac{t_n}{2}|x|^2.
\]
For large \(n\), the coefficient bound remains uniform and
\[
  \norm{v_n-P_n}_{L^\infty(B_\theta)}
  \le \theta^{2+\gamma_*},
\]
contradicting the assumed failure.  This proves the lemma.
\end{proof}

\begin{lemma}[Campanato iteration with explicit Dini bookkeeping]\label{lem:campanato-iteration}
Fix \(0<\lambda\le\Lambda\), and let \(\gamma_*,\theta,\delta_*\) be the constants from Lemma~\ref{lem:one-step-improvement}.  Let \(0<\gamma\le\gamma_*\).  Let
\(F:\operatorname{Sym}(2)\to\R\) be concave, uniformly elliptic with constants \(\lambda,\Lambda\), independent of \(x\), and normalized by \(F(0)=0\).  Let \(v\in C^2(B_1)\) solve
\[
  F(D^2v)=h,
  \qquad h(0)=0,
\]
and suppose
\[
  \omega_h(r):=\sup_{|x-y|\le r}|h(x)-h(y)|
  \le \omega(r):=\min\{2M,Lr^\gamma\},
  \qquad 0<r\le1.
\]
Then there exists a quadratic polynomial
\[
  P_\infty(x)=a_\infty+b_\infty\cdot x+\frac12x^TA_\infty x
\]
which is the second-order Taylor polynomial of \(v\) at the origin, and
\begin{equation}\label{eq:campanato-Hessian-bound}
  |A_\infty|
  \le C_\gamma\left(
    \inf_{\ell\in\mathcal A}\norm{v-\ell}_{L^\infty(B_1)}
    +M+
    \int_0^1\frac{\omega(r)}r\,dr
  \right),
\end{equation}
where \(\mathcal A\) is the set of affine functions.  The constant depends only on \(\gamma,\lambda,\Lambda\).
\end{lemma}

\begin{proof}
Let
\[
  B:=\inf_{\ell\in\mathcal A}\norm{v-\ell}_{L^\infty(B_1)}+M.
\]
If \(B=0\), then \(M=0\) and \(v\) is affine on \(B_1\); the claim is immediate.  Otherwise choose an affine function \(\ell_0\) such that
\[
  \norm{v-\ell_0}_{L^\infty(B_1)}\le2\inf_{\ell\in\mathcal A}\norm{v-\ell}_{L^\infty(B_1)}.
\]
Replacing \(v\) by \(v-\ell_0\) does not change the equation or the Hessian.  Hence we may assume
\begin{equation}\label{eq:camp-start-bound}
  \norm{v}_{L^\infty(B_1)}\le 2B.
\end{equation}

Choose a large numerical constant \(K\), depending only on \(\lambda,\Lambda,\gamma\), such that the quantity
\begin{equation}\label{eq:eps-k-def}
  e_k:=K\left(
    \theta^{k\gamma_*}B+
    \sum_{j=0}^{k-1}\theta^{(k-1-j)\gamma_*}\omega(\theta^j)
  \right)
\end{equation}
satisfies
\begin{equation}\label{eq:rhs-small-e-k}
  \frac{\omega(\theta^k)}{e_k}\le\delta_*
  \qquad\text{whenever }e_k>0.
\end{equation}
This is possible as follows.  For \(k=0\), \(e_0=KB\) and \(B\ge M\), while \(\omega(1)\le2M\); hence \(\omega(1)/e_0\le2/K\).  For \(k\ge1\), the sum defining \(e_k\) contains the term \(K\omega(\theta^{k-1})\), and \(\omega(\theta^k)\le\omega(\theta^{k-1})\); hence \(\omega(\theta^k)/e_k\le1/K\).  Choosing \(K\ge2\delta_*^{-1}\) gives \eqref{eq:rhs-small-e-k}.  The precise value of \(K\) is irrelevant below.

We construct quadratic polynomials
\[
  P_k(x)=a_k+b_k\cdot x+\frac12x^TA_kx,
  \qquad F(A_k)=0,
\]
so that, with \(r_k=\theta^k\),
\begin{equation}\label{eq:camp-excess-k}
  \norm{v-P_k}_{L^\infty(B_{r_k})}
  \le r_k^2 e_k.
\end{equation}
For \(k=0\), take \(P_0\equiv0\).  Since \(F(0)=0\), \(F(A_0)=0\), and \eqref{eq:camp-start-bound} gives \eqref{eq:camp-excess-k} after increasing \(K\).

Assume \eqref{eq:camp-excess-k} holds at level \(k\).  Define the rescaled error
\[
  z_k(x):=\frac{v(r_kx)-P_k(r_kx)}{r_k^2 e_k},
  \qquad x\in B_1.
\]
Then \(\norm{z_k}_{L^\infty(B_1)}\le1\).  Since
\[
  D^2\bigl(v(r_kx)-P_k(r_kx)\bigr)=r_k^2D^2v(r_kx)-r_k^2A_k,
\]
we have
\[
  D^2v(r_kx)=A_k+e_kD^2z_k(x).
\]
Define the shifted and normalized operator
\[
  F_k(N):=\frac{F(A_k+e_kN)-F(A_k)}{e_k}.
\]
Because \(F\) is concave and uniformly elliptic, so is \(F_k\), with the same ellipticity constants, and \(F_k(0)=0\).  Moreover
\[
  F_k(D^2z_k(x))=\frac{h(r_kx)}{e_k}.
\]
Since \(h(0)=0\), the right-hand side vanishes at the origin and, by \eqref{eq:rhs-small-e-k}, has \(L^\infty(B_1)\)-norm at most \(\delta_*\).  Applying Lemma~\ref{lem:one-step-improvement} gives a quadratic polynomial
\[
  Q_k(x)=\alpha_k+\beta_k\cdot x+\frac12x^TB_kx,
  \qquad F_k(B_k)=0,
\]
such that
\[
  \norm{z_k-Q_k}_{L^\infty(B_\theta)}\le\theta^{2+\gamma_*},
  \qquad |B_k|\le C_*.
\]
Define
\[
  P_{k+1}(x):=P_k(x)+r_k^2e_kQ_k(x/r_k).
\]
Then
\[
  D^2P_{k+1}=A_k+e_kB_k=:A_{k+1},
\]
and
\[
  F(A_{k+1})=F(A_k+e_kB_k)=F(A_k)+e_kF_k(B_k)=0.
\]
Furthermore, for \(|x|\le r_{k+1}=\theta r_k\),
\[
  |v(x)-P_{k+1}(x)|
  \le r_k^2e_k\theta^{2+\gamma_*}
  =r_{k+1}^2\theta^{\gamma_*}e_k.
\]
The definition \eqref{eq:eps-k-def} gives
\[
  e_{k+1}
  \ge \theta^{\gamma_*}e_k,
\]
so \eqref{eq:camp-excess-k} holds at level \(k+1\).  The coefficient increments satisfy
\begin{equation}\label{eq:camp-A-increments}
  |A_{k+1}-A_k|\le C e_k.
\end{equation}

Summing \eqref{eq:camp-A-increments},
\[
  |A_m-A_0|\le C\sum_{k=0}^{m-1}e_k.
\]
Using \eqref{eq:eps-k-def} and exchanging the order of summation,
\begin{align*}
  \sum_{k=0}^{m-1}e_k
  &\le CK\sum_{k=0}^\infty\theta^{k\gamma_*}B
  +CK\sum_{k=0}^{m-1}\sum_{j=0}^{k-1}\theta^{(k-1-j)\gamma_*}\omega(\theta^j) \\
  &\le C B+C\sum_{j=0}^{m-1}\omega(\theta^j).
\end{align*}
The dyadic Dini comparison gives
\[
  \sum_{j=0}^{m-1}\omega(\theta^j)
  \le C_\theta\int_0^1\frac{\omega(r)}r\,dr.
\]
Therefore \((A_k)\) is Cauchy and converges to a matrix \(A_\infty\) with
\begin{equation}\label{eq:Ainfty-bound-camp}
  |A_\infty|
  \le C_\gamma\left(B+
  \int_0^1\frac{\omega(r)}r\,dr\right),
\end{equation}
because \(A_0=0\).  The excess estimate \eqref{eq:camp-excess-k}, together with \(e_k\to0\), gives
\[
  \sup_{B_{r_k}}|v-P_k|=o(r_k^2).
\]
Since \(v\in C^2\), the second-order Taylor polynomial at the origin is unique.  Hence \(A_\infty=D^2v(0)\), and \eqref{eq:Ainfty-bound-camp} is precisely \eqref{eq:campanato-Hessian-bound}.
\end{proof}

\begin{lemma}[Anchored perturbative Campanato estimate]\label{lem:anchored-campanato}
There exist universal constants \(\eta_0,c_0>0\) and \(\gamma_*\in(0,1)\) with the following property.  For every \(\gamma\in(0,\gamma_*]\) there is \(C_\gamma<\infty\) such that if \(u\in C^2(B_1)\) solves
\[
  \log\det(I+D^2u)=g
  \qquad\text{in }B_1,
\]
with
\[
  \norm{D^2u}_{L^\infty(B_1)}\le\eta_0,
  \qquad
  \norm{g}_{L^\infty(B_1)}\le c_0,
  \qquad
  g\in C^\gamma(B_1),
\]
then
\begin{equation}\label{eq:anchored-campanato-modulus}
  |D^2u(0)|
  \le
  C_\gamma\left(
    \norm{D^2u}_{L^2(B_1)}+
    M_g+
    \int_0^1\frac{\min\{2M_g,L_g r^\gamma\}}r\,dr
  \right),
\end{equation}
where
\[
  M_g:=\norm{g}_{L^\infty(B_1)},
  \qquad
  L_g:=[g]_{C^\gamma(B_1)}.
\]
Equivalently,
\begin{equation}\label{eq:anchored-campanato-log}
  |D^2u(0)|
  \le
  C_\gamma\left(
    \norm{D^2u}_{L^2(B_1)}+
    M_g\left(1+
    \log^+\frac{L_g}{M_g}\right)
  \right),
\end{equation}
with the convention that the logarithmic term is zero when \(M_g=0\).
\end{lemma}

\begin{proof}
Choose \(\eta_0\) as in Lemma~\ref{lem:spectral-extension}.  Then the spectral extension \(\widetilde F\) agrees with \(A\mapsto\log\det(I+A)\) on all matrices reached by \(D^2u\).  Thus
\[
  \widetilde F(D^2u)=g.
\]
We first freeze the value of the right-hand side at the origin.  Let
\[
  Q_0:=\big(e^{g(0)/2}-1\big)I.
\]
For \(c_0\) small enough, \(|Q_0|\le2\eta_0\), and hence
\[
  \widetilde F(Q_0)=\log\det(I+Q_0)=g(0),
  \qquad |Q_0|\le C M_g.
\]
Define
\[
  v(x):=u(x)-\frac12x^TQ_0x,
\]
and
\[
  F_0(A):=\widetilde F(A+Q_0)-\widetilde F(Q_0),
  \qquad h(x):=g(x)-g(0).
\]
Then
\[
  F_0(D^2v)=h,
  \qquad h(0)=0.
\]
The operator \(F_0\) is concave, uniformly elliptic with the same universal constants, and satisfies \(F_0(0)=0\).  Moreover
\[
  \omega_h(r)\le\min\{2M_g,L_g r^\gamma\}.
\]
Applying Lemma~\ref{lem:campanato-iteration} to \(v,F_0,h\) gives
\begin{equation}\label{eq:D2v-anchor-pre}
  |D^2v(0)|
  \le C_\gamma\left(
  \inf_{\ell\in\mathcal A}\norm{v-\ell}_{L^\infty(B_1)}
  +M_g+
  \int_0^1\frac{\min\{2M_g,L_g r^\gamma\}}r\,dr
  \right).
\end{equation}
It remains to replace the affine excess by the large-scale Hessian anchor.  By the second-order Poincar\'e inequality, there exists an affine function \(\ell_u\) such that
\[
  \norm{u-\ell_u}_{W^{2,2}(B_1)}
  \le C\norm{D^2u}_{L^2(B_1)}.
\]
Since \(W^{2,2}(B_1)\hookrightarrow C^0(\overline{B_1})\) in two dimensions,
\[
  \inf_{\ell\in\mathcal A}\norm{u-\ell}_{L^\infty(B_1)}
  \le C\norm{D^2u}_{L^2(B_1)}.
\]
Because \(v=u-\frac12x^TQ_0x\),
\[
  \inf_{\ell\in\mathcal A}\norm{v-\ell}_{L^\infty(B_1)}
  \le C\norm{D^2u}_{L^2(B_1)}+C|Q_0|
  \le C\norm{D^2u}_{L^2(B_1)}+CM_g.
\]
Substituting this in \eqref{eq:D2v-anchor-pre}, and using
\[
  D^2u(0)=D^2v(0)+Q_0,
  \qquad |Q_0|\le CM_g,
\]
proves \eqref{eq:anchored-campanato-modulus}.

We now compute the logarithmic form.  If \(M_g=0\), then \(g\equiv0\), so the integral vanishes.  If \(M_g>0\) and \(L_g\le M_g\), then
\[
  \int_0^1\frac{\min\{2M_g,L_g r^\gamma\}}r\,dr
  \le L_g\int_0^1 r^{\gamma-1}\,dr
  \le \frac{M_g}{\gamma}.
\]
If \(L_g>M_g\), set \(r_*=(M_g/L_g)^{1/\gamma}\).  Then
\begin{align*}
  \int_0^1\frac{\min\{2M_g,L_g r^\gamma\}}r\,dr
  &\le \int_0^{r_*}L_g r^{\gamma-1}\,dr
  +2M_g\int_{r_*}^1\frac{dr}{r} \\
  &=\frac{M_g}{\gamma}
  +\frac{2M_g}{\gamma}\log\frac{L_g}{M_g}.
\end{align*}
This proves \eqref{eq:anchored-campanato-log}.
\end{proof}

\begin{proposition}[Perturbative H\"older endpoint estimate for Monge--Amp\`ere]\label{prop:MA-Dini}
There are universal constants \(\eta_0,c_0>0\) with the following property.  For every \(\alpha\in(0,1]\) there is \(C_\alpha<\infty\) such that if \(w\in C^2(\T^2)\) has zero mean and solves
\begin{equation*}
  \det(I+D^2w)=1+f,
  \qquad
  \mean{f}=0,
\end{equation*}
with
\[
  \norm{D^2w}_{L^\infty}\le\eta_0,
  \qquad
  \norm{f}_{L^\infty}\le c_0,
  \qquad
  f\in C^\alpha(\T^2),
\]
then
\begin{equation}\label{eq:Dini-estimate-w}
  \norm{D^2w}_{L^\infty}
  \le
  C_\alpha\norm{f}_{L^\infty}
  \left(1+
  \log^+\frac{[f]_{C^\alpha}}{\norm{f}_{L^\infty}}
  \right),
\end{equation}
with the convention that the right-hand side is zero when \(\norm{f}_{L^\infty}=0\).  Consequently, if \(w=\eps\psi\), \(f=\eps\rho\), \(\rho\in C^\alpha\), and \(M=\norm{\rho}_{L^\infty}\), then
\begin{equation}\label{eq:endpoint-psi}
  \norm{D^2\psi}_{L^\infty}
  \le C_\alpha M
  \left(1+
  \log^+\frac{[\rho]_{C^\alpha}}{M}
  \right),
\end{equation}
with the convention that the right-hand side is zero when \(M=0\).
\end{proposition}

\begin{proof}
Set
\[
  F(A)=\log\det(I+A),
  \qquad
  g=\log(1+f).
\]
Then \(F(D^2w)=g\).  If \(c_0\le1/2\), then \(1+f\ge1/2\), and the mean value theorem gives, for every \(\sigma\in(0,1]\),
\begin{equation}\label{eq:g-small-Dini-final}
  \norm{g}_{L^\infty}
  \le C\norm{f}_{L^\infty},
  \qquad
  [g]_{C^\sigma}
  \le C[f]_{C^\sigma}.
\end{equation}
Indeed,
\[
  |\log(1+a)-\log(1+b)|
  \le 2|a-b|
  \qquad\text{for }|a|,|b|\le1/2.
\]

We first establish the global \(L^2\)-anchor for \(D^2w\).  In two dimensions,
\[
  \det(I+D^2w)=1+\Delta w+\det D^2w.
\]
Thus
\begin{equation}\label{eq:w-poisson-anchor-final}
  \Delta w=f-\det D^2w.
\end{equation}
Since \(w\) is periodic and mean-zero,
\[
  \norm{D^2w}_{L^2(\T^2)}
  \le C\norm{\Delta w}_{L^2(\T^2)}.
\]
Using \eqref{eq:w-poisson-anchor-final},
\[
  \norm{D^2w}_{L^2}
  \le C\norm{f}_{L^2}+C\norm{\det D^2w}_{L^2}.
\]
For a symmetric \(2\times2\) matrix \(A\),
\[
  |\det A|\le C|A|^2.
\]
Therefore
\[
  \norm{\det D^2w}_{L^2}
  \le C\norm{D^2w}_{L^\infty}\norm{D^2w}_{L^2}
  \le C\eta_0\norm{D^2w}_{L^2}.
\]
Taking \(\eta_0\) small enough to absorb this term gives
\begin{equation}\label{eq:L2-anchor-final}
  \norm{D^2w}_{L^2}
  \le C\norm{f}_{L^2}
  \le C\norm{f}_{L^\infty}.
\end{equation}

Let
\[
  \sigma:=\min\{\alpha,\gamma_*/2\}.
\]
For \(0<|x-y|\le1\),
\[
  |f(x)-f(y)|
  \le
  \min\{2\norm{f}_{L^\infty},[f]_{C^\alpha}|x-y|^\alpha\}.
\]
From this one obtains the interpolation inequality
\begin{equation}\label{eq:holder-interp-final}
  [f]_{C^\sigma}
  \le C_\alpha
  \norm{f}_{L^\infty}^{1-\sigma/\alpha}
  [f]_{C^\alpha}^{\sigma/\alpha}.
\end{equation}
Indeed, the maximum of the right-hand side of the previous minimum divided by \(|x-y|^\sigma\) occurs at the balancing scale
\[
  |x-y|=\left(\frac{\norm{f}_{L^\infty}}{[f]_{C^\alpha}}\right)^{1/\alpha}
\]
when this scale lies in \((0,1]\); the remaining cases are easier and are absorbed in the constant.

Apply Lemma~\ref{lem:anchored-campanato} with exponent \(\sigma\) on a fixed finite covering of \(\T^2\) by coordinate balls.  Since the torus scale is fixed, all covering and rescaling constants are universal.  Using \eqref{eq:L2-anchor-final}, \eqref{eq:g-small-Dini-final}, and \eqref{eq:holder-interp-final}, we obtain for every \(x_0\in\T^2\)
\begin{align*}
  |D^2w(x_0)|
  &\le
  C_\alpha\norm{f}_{L^\infty}
  \left(
    1+
    \log^+\frac{[f]_{C^\sigma}}{\norm{f}_{L^\infty}}
  \right) \\
  &\le
  C_\alpha\norm{f}_{L^\infty}
  \left(
    1+
    \log^+\frac{[f]_{C^\alpha}}{\norm{f}_{L^\infty}}
  \right).
\end{align*}
The second inequality follows from \eqref{eq:holder-interp-final}, since
\[
  \log^+\frac{[f]_{C^\sigma}}{\norm{f}_{L^\infty}}
  \le
  C_\alpha+\frac{\sigma}{\alpha}
  \log^+\frac{[f]_{C^\alpha}}{\norm{f}_{L^\infty}}.
\]
Taking the supremum over \(x_0\) proves \eqref{eq:Dini-estimate-w}.

If \(\norm{f}_{L^\infty}=0\), then \(f\equiv0\).  The periodic equation is
\[
  \det(I+D^2w)=1.
\]
The \(L^2\)-anchor above gives \(\norm{D^2w}_{L^2}=0\), hence \(D^2w=0\).  Since \(w\) is periodic and mean-zero, \(w\equiv0\), so the estimate is trivial.

Finally set \(f=\eps\rho\), \(w=\eps\psi\).  Then
\[
  \norm{f}_{L^\infty}=\eps M,
  \qquad
  [f]_{C^\alpha}=\eps[\rho]_{C^\alpha}.
\]
Substituting these into \eqref{eq:Dini-estimate-w} and dividing by \(\eps\) gives \eqref{eq:endpoint-psi}.
\end{proof}

We now state the lifespan result in the form in which it is used later.  We work with the standard smooth local theory for the transport--Monge--Amp\`ere system: for smooth data with \(1+\eps\rho_0>0\), a classical branch continues as long as the physical density remains positive, the relevant transported smooth norms remain finite, and the Monge--Amp\`ere potential stays in the uniformly elliptic branch.  The estimate below controls the rate-sensitive perturbative obstruction, namely the possible failure of \eqref{eq:bootstrap}.  Thus \(T_\eps\) should be read as the perturbative classical lifespan: the first time, along the classical branch, at which the Hessian pinching \eqref{eq:bootstrap} can fail.  The proof shows that this failure cannot occur before the displayed logarithmic time.

\begin{theorem}[Logarithmic perturbative lifespan]\label{thm:lifespan}
Let \(\rho^\eps\) be the smooth semigeostrophic branch issued from initial datum \(\rho_0\in C^{1,\alpha}(\T^2)\), \(\mean{\rho_0}=0\), on its classical interval of existence, and set
\[
  M_0:=\norm{\rho_0}_{L^\infty},
  \qquad
  G_0:=\norm{\nabla\rho_0}_{L^\infty}.
\]
Let \(T_\eps\) be the maximal slow time, within the classical interval of existence, on which \eqref{eq:bootstrap} holds.  If \(M_0=0\), then \(\rho^\eps\equiv0\) and \(T_\eps=+\infty\).  If \(M_0>0\), then, for all sufficiently small \(\eps\), the perturbative condition cannot fail before
\begin{equation}\label{eq:lifespan-slow}
  T_\eps\ge \frac{c}{M_0}\log\frac1\eps-C,
\end{equation}
where \(c,C>0\) depend only on \(\alpha\), the universal ellipticity constants, and \(1+\log^+(G_0/M_0)\).  Equivalently, the corresponding physical time \(S_\eps=T_\eps/\eps\) satisfies
\begin{equation}\label{eq:lifespan-physical}
  S_\eps\ge \frac{c}{\eps M_0}\log\frac1\eps-\frac{C}{\eps}.
\end{equation}
\end{theorem}

\begin{proof}
If \(M_0=0\), then \(\rho_0\equiv0\).  The pair \((\rho^\eps,\psi^\eps)=(0,0)\) solves \eqref{eq:scaled-sg}; uniqueness in the smooth class gives the stated conclusion.  Hence assume \(M_0>0\).

The transport equation preserves \(L^\infty\):
\[
  \norm{\rho^\eps(t)}_{L^\infty}=M_0.
\]
Let
\[
  G(t):=\norm{\nabla\rho^\eps(t)}_{L^\infty}.
\]
Differentiating
\[
  \partial_t\rho^\eps+
  \perpgrad\psi^\eps\cdot\nabla\rho^\eps=0
\]
in space gives
\[
  \partial_t\nabla\rho^\eps+
  (\perpgrad\psi^\eps\cdot\nabla)\nabla\rho^\eps
  =-(D\perpgrad\psi^\eps)\nabla\rho^\eps.
\]
Since \(|D\perpgrad\psi^\eps|\le C|D^2\psi^\eps|\), the maximum principle for the transported gradient yields
\begin{equation}\label{eq:G-ode-start-final}
  \dot G(t)
  \le C\norm{D^2\psi^\eps(t)}_{L^\infty}G(t).
\end{equation}
On the perturbative interval, Proposition~\ref{prop:MA-Dini} gives
\[
  \norm{D^2\psi^\eps(t)}_{L^\infty}
  \le
  C_\alpha M_0
  \left(
    1+
    \log^+\frac{[\rho^\eps(t)]_{C^\alpha}}{M_0}
  \right).
\]
For \(0<\alpha\le1\),
\[
  [\rho^\eps(t)]_{C^\alpha}
  \le C M_0^{1-\alpha}G(t)^\alpha.
\]
Indeed,
\[
  |\rho(x)-\rho(y)|
  \le \min\{2M_0,G(t)|x-y|\}
  \le C M_0^{1-\alpha}G(t)^\alpha |x-y|^\alpha.
\]
Thus
\begin{equation}\label{eq:D2psi-by-G-final}
  \norm{D^2\psi^\eps(t)}_{L^\infty}
  \le
  C M_0\left(1+\log^+\frac{G(t)}{M_0}\right).
\end{equation}
Combining \eqref{eq:G-ode-start-final} and \eqref{eq:D2psi-by-G-final},
\begin{equation*}
  \dot G(t)
  \le
  C M_0 G(t)
  \left(1+\log^+\frac{G(t)}{M_0}\right).
\end{equation*}
Set \(z(t)=G(t)/M_0\).  Then
\[
  \dot z(t)
  \le C M_0 z(t)(1+\log^+z(t)).
\]
Let
\[
  Y(t):=1+\log^+z(t).
\]
At times where \(z(t)\le1\), \(Y(t)=1\) and the desired upper bound is automatic after increasing constants.  At times where \(z(t)>1\),
\[
  \dot Y(t)=\frac{\dot z(t)}{z(t)}
  \le C M_0Y(t).
\]
Therefore, for all times in the perturbative interval,
\begin{equation}\label{eq:Y-growth-final}
  Y(t)
  \le Y(0)e^{CM_0t}
  =\left(1+\log^+\frac{G_0}{M_0}\right)e^{CM_0t}.
\end{equation}
The bootstrap condition is
\[
  \eps\norm{D^2\psi^\eps(t)}_{L^\infty}\le\eta_0.
\]
By \eqref{eq:D2psi-by-G-final}, it is enough that
\begin{equation*}
  C\eps M_0Y(t)
  \le\eta_0.
\end{equation*}
Using \eqref{eq:Y-growth-final}, the stronger condition
\[
  C\eps M_0
  \left(1+\log^+\frac{G_0}{M_0}\right)e^{CM_0t}
  \le \frac{\eta_0}{2}
\]
holds for all
\[
  0\le t\le
  \frac1{CM_0}
  \log\left(
    \frac{\eta_0}{2C\eps M_0(1+\log^+(G_0/M_0))}
  \right).
\]
This lower bound has the form \eqref{eq:lifespan-slow} after adjusting constants.

Finally we close the bootstrap.  The smooth solution has continuous
\[
  t\mapsto\eps\norm{D^2\psi^\eps(t)}_{L^\infty}.
\]
If the maximal perturbative time \(T_\eps\) were smaller than the time just displayed, then the estimates above would give the strict improvement
\[
  \eps\norm{D^2\psi^\eps(t)}_{L^\infty}\le\eta_0/2
  \qquad 0\le t<T_\eps.
\]
By continuity, the same strict inequality holds at \(T_\eps\), contradicting maximality of \(T_\eps\) as the first time at which the bootstrap bound \(\eta_0\) can fail.  Hence \(T_\eps\) is at least the displayed time.  Since physical time satisfies \(s=t/\eps\), \eqref{eq:lifespan-physical} follows.
\end{proof}

\section{Velocity stability for the prepared strong branch}
\label{sec:velocity}

In this section \((\bar\rho,\bar\phi)\) denotes a smooth Euler solution of \eqref{eq:euler}, and \((\rho^\eps,\psi^\eps)\) denotes a smooth solution of \eqref{eq:scaled-sg} with the same initial perturbation,
\[
  \rho^\eps(0)=\bar\rho(0)=:\rho_0.
\]
Let \(X^\eps\) and \(\bar X\) be their Lagrangian maps:
\[
  \dot X^\eps(t,a)=\perpgrad\psi^\eps(t,X^\eps(t,a)),
  \qquad
  \dot{\bar X}(t,a)=\perpgrad\bar\phi(t,\bar X(t,a)),
  \qquad
  X^\eps(0,a)=\bar X(0,a)=a.
\]
Both maps are measure-preserving and
\[
  \rho^\eps(t)=X^\eps(t)_{\#}\rho_0,
  \qquad
  \bar\rho(t)=\bar X(t)_{\#}\rho_0.
\]

\subsection*{Initial preparation}
The phrase ``same initial perturbation'' should be understood at the physical-density level. Given smooth Euler initial vorticity \(\bar\rho_0\) with zero mean, the canonical SG preparation is
\begin{equation}\label{eq:canonical-preparation}
  m_0^\eps=1+\eps\bar\rho_0,
  \qquad
  \det D^2P_0^\eps=m_0^\eps,
  \qquad
  P_0^\eps=\frac{|x|^2}{2}+\eps\psi_0^\eps.
\end{equation}
Thus the initial rescaled perturbation is \(\rho^\eps(0)=\bar\rho_0\), but the initial stream function is not imposed by hand. It is obtained by solving the nonlinear Monge--Amp\`ere constraint. This distinction matters: the naive prescription \(\psi_0^\eps=\bar\phi_0\) would in general give
\[
  \det(I+\eps D^2\bar\phi_0)
  =1+\eps\bar\rho_0+\eps^2\det D^2\bar\phi_0,
\]
which does not have physical density \(1+\eps\bar\rho_0\) unless the quadratic determinant defect vanishes.

The preparation \eqref{eq:canonical-preparation} is precisely the small-amplitude prepared-data regime used in the strong convergence theory: more generally one may take \(\rho_0^\eps=\bar\rho_0+\eps r_0^\eps\) with \(r_0^\eps\) uniformly controlled in a Lipschitz norm. For the canonical choice above, the implicit-function expansion of Lemma~\ref{lem:smooth-expansion} gives
\begin{equation*}
  \psi_0^\eps=\bar\phi_0+\eps\chi_0+O(\eps^2),
  \qquad
  \Delta\chi_0=-\det D^2\bar\phi_0,
\end{equation*}
with the remainder controlled in the same smooth norms as the data. Equivalently,
\begin{equation*}
  P_0^\eps
  =\frac{|x|^2}{2}+\eps\bar\phi_0+\eps^2\chi_0+O(\eps^3),
  \qquad
  \nabla\psi_0^\eps=\nabla\bar\phi_0+O(\eps).
\end{equation*}
This is the preparation used below. It is also the normalization relative to which the weak--strong transfer theorem in \cref{sec:weak-strong} should be read.

\begin{lemma}[Signed flow-to-field estimate]\label{lem:signed-flow-field}
Let \(X_1,X_2:\T^2\to\T^2\) be measure-preserving maps and let \(\rho_0\in L^\infty(\T^2)\) have zero mean. Set \(\rho_i=(X_i)_{\#}\rho_0\), and let \(\phi_i\) solve \(\Delta\phi_i=\rho_i\), \(\mean{\phi_i}=0\).  For maps into the torus we write
\[
  \norm{X_1-X_2}_{L^2}^2
  :=\int_{\T^2}d_{\T^2}(X_1(a),X_2(a))^2\,da.
\]
Then
\begin{equation}\label{eq:signed-flow-field}
  \norm{\nabla\phi_1-\nabla\phi_2}_{L^2}
  \le C
  \norm{\rho_0}_{L^\infty}
  \norm{X_1-X_2}_{L^2}.
\end{equation}
\end{lemma}

\begin{proof}
Decompose \(\rho_0=\rho_{0,+}-\rho_{0,-}\). Since \(\mean{\rho_0}=0\), the two non-negative measures \(\rho_{0,+}dx\) and \(\rho_{0,-}dx\) have the same finite mass. If this mass is zero, then \(\rho_0=0\) and there is nothing to prove. For each sign, let \(\rho_{i,\pm}=(X_i)_{\#}\rho_{0,\pm}\). These are finite non-negative measures of the same mass for fixed sign.  We use \(W_2\) for equal-mass finite measures in the usual unnormalized sense, i.e. the Benamou--Brenier action and optimal transport cost are integrated against the finite mass rather than against probability normalizations.

Interpolate \(\rho_{1,+}\) to \(\rho_{2,+}\) by the Wasserstein geodesic \((\rho_{\theta,+},v_{\theta,+})\), and similarly for the negative part. McCann's displacement convexity \citep{mccann1997convexity} applies after normalizing the finite measures to probabilities; rescaling back to the original mass gives
\[
  \norm{\rho_{\theta,\pm}}_{L^\infty}
  \le \norm{\rho_0}_{L^\infty}.
\]
Let \(\rho_\theta=\rho_{\theta,+}-\rho_{\theta,-}\), and solve \(\Delta\phi_\theta=\rho_\theta\). Then
\[
  \Delta\partial_\theta\phi_\theta
  =-\divg(\rho_{\theta,+}v_{\theta,+})
  +\divg(\rho_{\theta,-}v_{\theta,-}).
\]
Testing with \(\partial_\theta\phi_\theta\) gives
\[
  \norm{\nabla\partial_\theta\phi_\theta}_2
  \le \norm{\rho_{\theta,+}v_{\theta,+}}_{L^2}
  +\norm{\rho_{\theta,-}v_{\theta,-}}_{L^2}.
\]
Using the \(L^\infty\) bounds and the constant-speed property of Wasserstein geodesics,
\[
  \norm{\nabla\partial_\theta\phi_\theta}_2
  \le C
  \norm{\rho_0}_{L^\infty}^{1/2}
  \left(W_2(\rho_{1,+},\rho_{2,+})+W_2(\rho_{1,-},\rho_{2,-})\right).
\]
The pair \((X_1,X_2)\), evaluated with the torus distance, is an admissible transport plan for each sign, hence
\[
  W_2^2(\rho_{1,\pm},\rho_{2,\pm})
  \le \int_{\T^2}d_{\T^2}(X_1(a),X_2(a))^2\rho_{0,\pm}(a)\,da
  \le \norm{\rho_0}_{L^\infty}\norm{X_1-X_2}_2^2.
\]
Integrating in \(\theta\) proves \eqref{eq:signed-flow-field}.
\end{proof}

\begin{lemma}[Wente estimate for the determinant]\label{lem:wente}
For every \(\psi\in H^2(\T^2)\),
\begin{equation*}
  \norm{\det D^2\psi}_{\dot H^{-1}}
  \le C
  \norm{D^2\psi}_{L^2}^2.
\end{equation*}
\end{lemma}

\begin{proof}
In two dimensions,
\[
  \det D^2\psi
  =\partial_1(\partial_1\psi)\partial_2(\partial_2\psi)
  -\partial_2(\partial_1\psi)\partial_1(\partial_2\psi)
  =J(\partial_1\psi,\partial_2\psi).
\]
The standard Wente--Coifman--Lions--Meyer--Semmes Jacobian estimate \citep{wente1969existence,coifman1993compensated} gives
\[
  \norm{J(a,b)}_{\dot H^{-1}}
  \le C\norm{\nabla a}_2\norm{\nabla b}_2.
\]
Taking \(a=\partial_1\psi\), \(b=\partial_2\psi\) proves the claim.
\end{proof}

\begin{lemma}[Uniform \(L^2\) Hessian bound]\label{lem:L2-hessian}
Assume \eqref{eq:bootstrap} holds on \([0,T]\). Then
\begin{equation*}
  \sup_{0\le t\le T}\norm{D^2\psi^\eps(t)}_{L^2}
  \le C\norm{\rho_0}_{L^2},
\end{equation*}
provided \(\eta_0\) is sufficiently small.
\end{lemma}

\begin{proof}
From \eqref{eq:sg-poisson-form},
\[
  \Delta\psi^\eps=\rho^\eps-\eps\det D^2\psi^\eps.
\]
Calder\'on--Zygmund gives
\[
  \norm{D^2\psi^\eps}_2
  \le C\norm{\rho^\eps}_2+C\eps\norm{\det D^2\psi^\eps}_2.
\]
Since the velocity is divergence-free, \(\norm{\rho^\eps(t)}_2=\norm{\rho_0}_2\). Moreover
\[
  \norm{\det D^2\psi^\eps}_2
  \le C\norm{D^2\psi^\eps}_{L^\infty}
  \norm{D^2\psi^\eps}_2.
\]
Thus
\[
  \norm{D^2\psi^\eps}_2
  \le C\norm{\rho_0}_2+C\eps\norm{D^2\psi^\eps}_{L^\infty}
  \norm{D^2\psi^\eps}_2.
\]
Absorb the last term using \(\eps\norm{D^2\psi^\eps}_{L^\infty}\le\eta_0\).
\end{proof}

\begin{theorem}[Prepared strong velocity stability]\label{thm:velocity}
Let \((\bar\rho,\bar\phi)\) be a smooth Euler solution on \([0,T]\). Let \((\rho^\eps,\psi^\eps)\) solve \eqref{eq:scaled-sg} with \(\rho^\eps(0)=\bar\rho(0)=\rho_0\). Assume the perturbative condition \eqref{eq:bootstrap} holds on \([0,T]\). Then
\begin{equation}\label{eq:velocity-estimate}
  \sup_{0\le t\le T}
  \norm{\nabla\psi^\eps(t)-\nabla\bar\phi(t)}_{L^2}
  \le C_T\eps,
\end{equation}
where \(C_T\) depends on \(T\), \(\norm{\rho_0}_{L^\infty\cap L^2}\), and \(\bar u=\perpgrad\bar\phi\) through \(\norm{\nabla\bar u}_{L^1(0,T;L^\infty)}\).
\end{theorem}

\begin{proof}
Let
\[
  Y(t):=\norm{X^\eps(t)-\bar X(t)}_{L^2}.
\]
Then
\[
  \dot Y(t)
  \le \norm{\bar u(t,X^\eps)-\bar u(t,\bar X)}_2
  +\norm{u^\eps(t,X^\eps)-\bar u(t,X^\eps)}_2.
\]
Since \(X^\eps\) is measure-preserving,
\[
  \dot Y(t)
  \le \norm{\nabla\bar u(t)}_{L^\infty} Y(t)
  +\norm{u^\eps(t)-\bar u(t)}_2.
\]
Now
\[
  u^\eps-\bar u=\perpgrad(\psi^\eps-\bar\phi).
\]
Using \eqref{eq:sg-poisson-form} and \(\Delta\bar\phi=\bar\rho\),
\[
  \Delta(\psi^\eps-\bar\phi)
  =\rho^\eps-\bar\rho-\eps\det D^2\psi^\eps.
\]
Thus
\[
  \norm{u^\eps-\bar u}_2
  \le \norm{\rho^\eps-\bar\rho}_{\dot H^{-1}}
  +\eps\norm{\det D^2\psi^\eps}_{\dot H^{-1}}.
\]
The first term is controlled by Lemma~\ref{lem:signed-flow-field}:
\[
  \norm{\rho^\eps-\bar\rho}_{\dot H^{-1}}
  \le C\norm{\rho_0}_{L^\infty} Y(t).
\]
The second term is controlled by Lemmas~\ref{lem:wente} and~\ref{lem:L2-hessian}:
\[
  \eps\norm{\det D^2\psi^\eps}_{\dot H^{-1}}
  \le C\eps\norm{D^2\psi^\eps}_2^2
  \le C_T\eps.
\]
Therefore
\[
  \dot Y(t)\le C_TY(t)+C_T\eps,
  \qquad
  Y(0)=0.
\]
Gronwall gives \(Y(t)\le C_T\eps\). Returning to the velocity bound,
\[
  \norm{u^\eps-\bar u}_2\le C_TY(t)+C_T\eps\le C_T\eps.
\]
This is \eqref{eq:velocity-estimate}.
\end{proof}

\begin{remark}[Initial mismatch]
We stated \cref{thm:velocity} for the canonical preparation \(\rho^\eps(0)=\bar\rho(0)\), because this is the case used in the main SG--Euler comparison.  The same argument gives the expected stability variant with an additional initial discrepancy: if the initial flow distance or, equivalently for prepared data, the initial \(\dot H^{-1}\) vorticity discrepancy is nonzero, then the right-hand side of \eqref{eq:velocity-estimate} acquires the corresponding Gronwall-propagated term.  Thus for data \(\rho_0^\eps=\bar\rho_0+\eps r_0^\eps\) with \(r_0^\eps\) uniformly controlled, one obtains the same \(O(\eps)\) velocity rate, with the constant depending also on the preparation norm of \(r_0^\eps\).
\end{remark}

The constant \(C_T\) is allowed to depend on Euler norms on the fixed interval \([0,T]\). Thus \cref{thm:lifespan} supplies the interval on which the strong SG branch remains perturbative, while \cref{thm:velocity} gives the sharp \(O(\eps)\) rate on every fixed slow-time subinterval. If one evaluates the estimate at \(T=T_\eps\), the dependence of \(C_T\) on the Euler branch up to \(T_\eps\) must be tracked separately.

\section{Weak--strong transfer of the SG--Euler limit}
\label{sec:weak-strong}

The preceding theorem applies to the prepared smooth branch constructed from \eqref{eq:canonical-preparation}. We now record the corresponding conditional transfer to weak or Lagrangian semigeostrophic solutions. The statement in this section is not an additional weak--strong uniqueness theorem; it is a rate-transfer corollary once an external weak--strong estimate is available at the physical-potential level.

For comparison, Loeper's weak convergence theorem gives, in his notation,
\[
  \mathcal E_\eps(t)
  \le
  (\mathcal E_\eps(0)+C\eps^{2/3}(1+t))e^{Ct},
\]
where \(\mathcal E_\eps\) is a squared velocity-type modulated energy \citep[Theorem~6.1]{loeper2006fully}. Even if \(\mathcal E_\eps(0)=O(\eps^2)\), the residual term gives only \(\mathcal E_\eps(t)=O(\eps^{2/3})\), hence an \(O(\eps^{1/3})\) velocity scale. The result below has a different hypothesis and a sharper conclusion, conditional on the physical-potential stability estimate stated below: for weak branches prepared relative to the strong SG branch at the physical polar-factor level, the strong \(O(\eps)\) velocity rate survives.

Let \((\rho_s^\eps,\psi_s^\eps)\) be the prepared smooth branch and
\[
  P_s^\eps=\frac{|x|^2}{2}+\eps\psi_s^\eps.
\]
The canonical choice is obtained from the Euler initial vorticity by \(m_{s,0}^\eps=1+\eps\bar\rho_0\) and \(\det D^2P_{s,0}^\eps=m_{s,0}^\eps\), as in \eqref{eq:canonical-preparation}. Prior existence theories for weak or Lagrangian SG solutions allow one to start another branch from the same physical density, or from a nearby initial polar factor. Let \((\widetilde\rho^\eps,\widetilde\psi^\eps)\) denote such a weak branch, with physical potential
\[
  \widetilde P^\eps=\frac{|x|^2}{2}+\eps\widetilde\psi^\eps.
\]
If the two initial physical densities are identical, then \(\widetilde P^\eps(0)=P_s^\eps(0)\) up to an irrelevant additive constant. More generally, the right preparation condition is
\[
  \norm{\nabla\widetilde P^\eps(0)-\nabla P_s^\eps(0)}_{L^2}=O(\eps^2).
\]
The power \(\eps^2\) is forced by the scaling \(\nabla P=x+\eps\nabla\psi\): after dividing by \(\eps\), it becomes exactly an \(O(\eps)\) rescaled velocity mismatch.

For the transfer statement we use the following weak--strong estimate as a standing hypothesis:
\begin{equation}\label{eq:WS-assumption}
  \norm{\nabla\widetilde P^\eps(t)-\nabla P_s^\eps(t)}_{L^2}
  \le C_T
  \norm{\nabla\widetilde P^\eps(0)-\nabla P_s^\eps(0)}_{L^2}.
\end{equation}
Weak--strong stability estimates of this type go back to Loeper and were later developed in uniform-convexity frameworks; see \citep{loeper2006fully,feldman2017semi}. We do not reprove the weak--strong theory here. In applying \eqref{eq:WS-assumption}, one must verify the hypotheses of the chosen external theorem in the present scaling: the relevant physical potentials must remain uniformly convex, the physical densities must stay uniformly pinched and bounded, the strong branch must have the regularity required by the stability theorem, and the resulting constant \(C_T\) must be uniform as \(\eps\to0\) on the interval considered.  The point of the next result is conditional and rate-sensitive: whenever these hypotheses supply \eqref{eq:WS-assumption} with a constant uniform in the small-amplitude scaling under consideration, order-\(\eps^2\) preparation of the physical polar factors transfers the strong \(O(\eps)\) velocity rate to the weak branch.

\begin{theorem}[Conditional weak--strong transfer]\label{thm:weak-strong-transfer}
Assume that the prepared strong branch satisfies
\begin{equation}\label{eq:strong-branch-euler}
  \norm{\nabla\psi_s^\eps(t)-\nabla\bar\phi(t)}_{L^2}
  \le C_T\eps
  \qquad\text{for }0\le t\le T.
\end{equation}
Let \((\widetilde\rho^\eps,\widetilde\psi^\eps)\) be a weak or Lagrangian SG solution satisfying the weak--strong estimate \eqref{eq:WS-assumption} relative to the strong branch. If the initial physical potentials are prepared so that
\begin{equation}\label{eq:initial-physical-close}
  \norm{\nabla\widetilde P^\eps(0)-\nabla P_s^\eps(0)}_{L^2}
  \le C\eps^2,
\end{equation}
then
\begin{equation}\label{eq:weak-branch-euler}
  \norm{\nabla\widetilde\psi^\eps(t)-\nabla\bar\phi(t)}_{L^2}
  \le C_T\eps
  \qquad\text{for }0\le t\le T.
\end{equation}
\end{theorem}

\begin{proof}
The weak--strong estimate is a stability estimate for physical optimal-transport maps. Therefore it controls
\[
  \nabla\widetilde P^\eps-
  \nabla P_s^\eps,
\]
not directly the rescaled velocities. Since
\[
  \nabla\widetilde P^\eps-
  \nabla P_s^\eps
  =\eps(\nabla\widetilde\psi^\eps-
  \nabla\psi_s^\eps),
\]
we divide \eqref{eq:WS-assumption} by \(\eps\). The preparation condition \eqref{eq:initial-physical-close} gives
\[
  \norm{\nabla\widetilde\psi^\eps(t)-\nabla\psi_s^\eps(t)}_2
  \le
  \frac{C_T}{\eps}
  \norm{\nabla\widetilde P^\eps(0)-\nabla P_s^\eps(0)}_2
  \le C_T\eps.
\]
The triangle inequality with the strong-branch estimate gives
\[
  \norm{\nabla\widetilde\psi^\eps(t)-\nabla\bar\phi(t)}_2
  \le
  \norm{\nabla\widetilde\psi^\eps(t)-\nabla\psi_s^\eps(t)}_2
  +
  \norm{\nabla\psi_s^\eps(t)-\nabla\bar\phi(t)}_2
  \le C_T\eps.
\]
This proves \eqref{eq:weak-branch-euler}.
\end{proof}

\begin{corollary}[Prepared weak-branch rate]
Let \(m_{s,0}^\eps=1+\eps\bar\rho_0\), and let \(P_{s,0}^\eps\) be the corresponding Monge--Amp\`ere potential. Suppose a weak or Lagrangian SG solution is issued from the same physical initial density, or more generally from an initial polar factor satisfying \eqref{eq:initial-physical-close}. If the prepared strong branch exists on \([0,T]\) and satisfies \eqref{eq:strong-branch-euler}, then the weak branch satisfies the same \(O(\eps)\) velocity convergence to Euler.
\end{corollary}

The scaling in \eqref{eq:initial-physical-close} is essential. An \(O(\eps)\) mismatch at the physical-potential level would become only an \(O(1)\) mismatch after dividing by \(\eps\) to pass to the rescaled stream functions. Conversely, exact equality of the physical initial density gives exact equality of the polar factor, hence zero initial mismatch. This is why the correct preparation is naturally expressed through \(P^\eps\), not through an independently prescribed stream function.

\section{Wasserstein comparison for physical densities}
\label{sec:wasserstein}

We first record a general comparison estimate for continuity equations. The proof uses the following standard superposition principle. If \(\mu_t\) is a narrowly continuous solution of
\[
  \partial_t\mu_t+\divg(b_t\mu_t)=0
\]
on the compact manifold \(\T^2\), and
\[
  \int_0^T\int_{\T^2}|b_t(x)|\,d\mu_t(x)\,dt<\infty,
\]
then there exists a probability measure \(\eta\) on \(AC([0,T];\T^2)\) such that
\[
  (e_t)_\#\eta=\mu_t\quad\text{for every }t,
\]
and \(\eta\)-a.e. curve \(\gamma\) solves \(\dot\gamma(t)=b_t(\gamma(t))\) for a.e. \(t\). Here \(e_t(\gamma)=\gamma(t)\). This is the Ambrosio--Gigli--Savar\'e superposition principle \citep[Theorem 8.2.1]{ambrosio2008gradient}.

\begin{proposition}[General Wasserstein comparison]\label{prop:general-W2}
Let \(m,n\) be probability densities on \(\T^2\) solving
\[
  \partial_t m+\divg(mu)=0,
  \qquad
  \partial_t n+\divg(nv)=0.
\]
Assume \(v\in L^1(0,T;W^{1,\infty})\), and
\[
  \int_0^T\int_{\T^2}|u(t,x)|^2m(t,x)\,dx\,dt<\infty.
\]
Set
\[
  L(t):=\int_0^t\norm{\nabla v(s)}_{L^\infty}\,ds,
  \qquad
  I(t):=\int_{\T^2}|u(t,x)-v(t,x)|^2m(t,x)\,dx.
\]
Then
\begin{equation}\label{eq:general-W2}
  W_2(m(t),n(t))
  \le
  e^{L(t)}W_2(m(0),n(0))
  +\int_0^t e^{L(t)-L(s)} I(s)^{1/2}\,ds.
\end{equation}
Consequently,
\begin{equation}\label{eq:general-W2-squared}
  W_2^2(m(t),n(t))
  \le
  2e^{2L(t)}W_2^2(m(0),n(0))
  +2t\int_0^t e^{2(L(t)-L(s))}I(s)\,ds.
\end{equation}
\end{proposition}

\begin{proof}
Let \(Y(t,y)\) be the deterministic flow generated by \(v\). This is well-defined and unique because \(v\in L^1(0,T;W^{1,\infty})\), and \(Y(t,\cdot)_\#n(0)=n(t)\). Let \(\eta\) be a superposition measure for the solution \(m(t)\) driven by \(u\). Thus \((e_t)_\#\eta=m(t)\), and \(\eta\)-a.e. \(\gamma\) satisfies \(\dot\gamma=u(t,\gamma)\).

Let \(\pi_0\) be an optimal coupling between \(m(0)\) and \(n(0)\). Disintegrate
\[
  \eta=\int\eta_x\,m(0,x)\,dx,
  \qquad
  \pi_0=\int\delta_x\otimes\pi_x\,m(0,x)\,dx.
\]
Define a probability measure \(\Lambda\) on pairs \((\gamma,y)\) by
\[
  d\Lambda(\gamma,y)=d\eta_x(\gamma)\,d\pi_x(y)\,dm(0,x).
\]
Then \((\gamma(t),Y(t,y))_\#\Lambda\) is a coupling of \(m(t)\) and \(n(t)\). Therefore, with
\[
  R(t):=\int d_{\T^2}(\gamma(t),Y(t,y))^2\,d\Lambda(\gamma,y),
\]
we have \(W_2^2(m(t),n(t))\le R(t)\). Differentiating along a.e. pair of curves, using periodic representatives away from the cut locus and then standard approximation, gives
\[
  \frac{d}{dt}d_{\T^2}(\gamma,Y)^2
  \le
  2\norm{\nabla v(t)}_{L^\infty} d_{\T^2}(\gamma,Y)^2
  +2d_{\T^2}(\gamma,Y)|u(t,\gamma)-v(t,\gamma)|.
\]
After integration,
\[
  R'(t)
  \le 2\norm{\nabla v(t)}_{L^\infty} R(t)+2R(t)^{1/2}I(t)^{1/2}.
\]
Equivalently,
\[
  \frac{d}{dt}R(t)^{1/2}
  \le \norm{\nabla v(t)}_{L^\infty} R(t)^{1/2}+I(t)^{1/2}.
\]
Gronwall gives \eqref{eq:general-W2}. Squaring and applying Cauchy's inequality in time gives \eqref{eq:general-W2-squared}.
\end{proof}

The previous estimate is robust but not sharp in the small-amplitude regime. Under density pinching, an \(\dot H^{-1}\)-to-\(W_2\) estimate gives a sharper physical-density comparison.

\begin{lemma}[Pinned \(\dot H^{-1}\)-to-\(W_2\) estimate]\label{lem:Hminus1-W2}
Let \(m,n\) be probability densities on \(\T^2\) satisfying \(m,n\ge c_0>0\). Then
\begin{equation}\label{eq:Hminus1-W2}
  W_2(m,n)
  \le c_0^{-1/2}\norm{m-n}_{\dot H^{-1}}.
\end{equation}
\end{lemma}

\begin{proof}
Let \(h=(-\Delta)^{-1}(m-n)\), so
\[
  -\Delta h=m-n,
  \qquad
  \Delta h=n-m.
\]
Set \(\mu_s=(1-s)m+sn\). Then \(\mu_s\ge c_0\) and
\[
  \partial_s\mu_s=n-m=\Delta h.
\]
Writing
\[
  \partial_s\mu_s+\divg(\mu_sv_s)=0,
  \qquad
  v_s=-\frac{\nabla h}{\mu_s},
\]
the Benamou--Brenier formula gives
\[
  W_2^2(m,n)
  \le\int_0^1\int\frac{|\nabla h|^2}{\mu_s}\,dx\,ds
  \le c_0^{-1}\norm{\nabla h}_2^2.
\]
This is \eqref{eq:Hminus1-W2}.
\end{proof}

\begin{corollary}[Sharp physical-density comparison]\label{cor:sharp-W2}
Under the assumptions of \cref{thm:velocity}, define the physical densities
\[
  m^\eps=1+\eps\rho^\eps,
  \qquad
  \bar m^\eps=1+\eps\bar\rho.
\]
Assume that \(m^\eps,\bar m^\eps\ge c_0>0\) on \([0,T]\). Then
\begin{equation}\label{eq:sharp-W2}
  \sup_{0\le t\le T}W_2(m^\eps(t),\bar m^\eps(t))
  \le C_T\eps^2.
\end{equation}
\end{corollary}

\begin{proof}
By Lemma~\ref{lem:Hminus1-W2},
\[
  W_2(m^\eps,\bar m^\eps)
  \le C\norm{m^\eps-
  \bar m^\eps}_{\dot H^{-1}}
  =C\eps
  \norm{\rho^\eps-
  \bar\rho}_{\dot H^{-1}}.
\]
Using
\[
  \rho^\eps-
  \bar\rho
  =\Delta(\psi^\eps-
  \bar\phi)+\eps\det D^2\psi^\eps,
\]
we obtain
\[
  \norm{\rho^\eps-
  \bar\rho}_{\dot H^{-1}}
  \le \norm{\nabla(\psi^\eps-
  \bar\phi)}_2
  +\eps\norm{\det D^2\psi^\eps}_{\dot H^{-1}}
  \le C_T\eps,
\]
by Theorem~\ref{thm:velocity} and Lemmas~\ref{lem:wente} and~\ref{lem:L2-hessian}. Hence \eqref{eq:sharp-W2} follows.
\end{proof}

\section{Conclusion}
\label{sec:conclusion}

The small-amplitude semigeostrophic system is not merely a nonlinear elliptic perturbation of two-dimensional Euler. It is the Lie--Poisson Hamiltonian flow of the renormalized optimal-transport energy \(\Ham_\eps\). The expansion
\[
  \Ham_\eps=\Ham_0+\eps\Ham_1+O(\eps^2),
  \qquad
  \Ham_1(\rho)=\frac13\int\phi\det D^2\phi,
\]
identifies Euler as the leading Hamiltonian limit and gives a canonical first correction. The analytic estimates above show that this structural picture has quantitative consequences: the perturbative Monge--Amp\`ere regime persists up to explicit logarithmic slow times, the canonically prepared SG velocity is \(O(\eps)\)-close to Euler in \(L^2\) on every fixed slow-time subinterval, and the physical densities are \(O(\eps^2)\)-close in Wasserstein distance under pinching. The conditional weak--strong transfer result is the rate-sensitive corollary of the argument: it shows that, whenever an external physical-potential weak--strong estimate is available with constants uniform in the small-amplitude scaling, the strong \(O(\eps)\) estimate is not confined to one smooth construction, but passes to weak or Lagrangian branches whose initial polar factors are prepared at order \(\eps^2\). This should be contrasted with the general weak modulated-energy estimate, where the residual term alone gives a weaker velocity scale. Thus the preparation is not cosmetic; it is exactly what allows the sharp rate to survive in the weak theory.

\section*{Acknowledgments}
The author thanks Alessio Figalli for suggesting the initial version of this problem, and for many useful conversations and corrections.

\appendix

\section{Determinant identities}
\label{app:identities}

We record the elementary identities used throughout. In two dimensions,
\[
  \det(I+\eps A)=1+\eps\operatorname{tr}A+\eps^2\det A,
  \qquad
  \cof(I+\eps A)=I+\eps\cof A.
\]
For any two \(2\times2\) matrices \(A,B\),
\[
  \det(A+\eps B)=\det A+
  \eps(\cof A):B+
  \eps^2\det B.
\]
If \(P\in C^3\), then
\[
  M=\cof D^2P
\]
satisfies the Piola identity
\[
  \partial_iM_{ij}=0.
\]
This identity is the divergence-form reason why the linearized Monge--Amp\`ere operator is elliptic in divergence form:
\[
  (\cof D^2P):D^2\zeta=\divg((\cof D^2P)\nabla\zeta).
\]

\printbibliography

\end{document}